\documentclass{ifacconf}

\usepackage{natbib}        % required for bibliography
\usepackage[left=0.6in,right=0.5565in,bottom=0.55in,top=1.45in]{geometry}

\usepackage{amsmath,amssymb,amsfonts}
\usepackage{graphicx}
\usepackage{textcomp}
\usepackage{xcolor}
\usepackage{booktabs}
 \usepackage{psfrag}
\usepackage[width=.75\textwidth]{caption}
\usepackage{subcaption}
\usepackage{balance}
\usepackage{multirow}
\usepackage{tabularx}

\usepackage[linesnumbered,ruled,vlined,algo2e]{algorithm2e}

\SetKwInput{KwInput}{Input}                % Set the Input
\SetKwInput{KwOutput}{Output}              % set the Output

\usepackage{dsfont}
\usepackage{tikz}
\usetikzlibrary{shapes.geometric}
\usepackage{pgfplots}
\pgfplotsset{width=0.45\textwidth}

\usepackage{graphicx}
\graphicspath{{Figures/}}

\newcommand{\T}{^\mathsf{T}} % matrix transpose
\renewcommand{\d}{\mathrm{d}} % derivative operator
\newcommand{\R}{\mathds{R}} % real numbers
\newcommand{\N}{\mathds{N}} % natural numbers
\DeclareMathOperator*{\argmax}{arg\,max}

\usepackage[noend]{algpseudocode}
\usepackage[nothing]{algorithm}

\algnewcommand\And{\textbf{and}}

\usepackage{array}
\newcolumntype{P}[1]{>{\centering\arraybackslash}p{#1}}

\begin{document}
\begin{frontmatter}

\title{Efficient sample selection for safe learning\thanksref{footnoteinfo}}

\thanks[footnoteinfo]{Research supported by NCCR Automation, a National Centre of Competence in Research, funded by the Swiss National Science Foundation (grant no. 180545), and by the European Research Council (ERC) under the H2020 Advanced Grant no. 787845 (OCAL).}

\author[ETH]{Marta Zagorowska} 
\author[ETH]{Efe C. Balta} 
\author[ETH,EMPA]{Varsha Behrunani} 
\author[ETH]{Alisa Rupenyan} 
\author[ETH]{John Lygeros} 

\address[ETH]{Automatic Control Laboratory, ETH Zurich, Switzerland, (E-mails: \{mzagorowska,ebalta,bvarsha,ralisa,lygeros\}@control.ee.ethz.ch)}
\address[EMPA]{Urban Energy Systems Laboratory, Swiss Federal Laboratories for Materials Science and
Technology (EMPA), Dübendorf, Switzerland}

\begin{abstract}                % Abstract of not more than 250 words.
Ensuring safety in industrial control systems usually involves imposing constraints at the design stage of the control algorithm. Enforcing constraints is challenging if the underlying functional form is unknown. The challenge can be addressed by using surrogate models, such as Gaussian processes, which provide confidence intervals used to find solutions that can be considered safe. This in turn involves an exhaustive search on the entire search space. That approach can quickly become computationally expensive. We reformulate the exhaustive search as a series of optimization problems to find the next recommended points. We show that the proposed reformulation allows using a wide range of available optimization solvers, such as derivative-free methods. We show that by exploiting the properties of the solver, we enable the introduction of new stopping criteria into safe learning methods and increase flexibility in trading off solver accuracy and computational time. The results from a non-convex optimization problem and an application for controller tuning confirm the flexibility and the performance of the proposed reformulation.

\end{abstract}


\end{frontmatter}
%===============================================================================

\section{Introduction}
Ensuring safety in industrial control systems is usually done by enforcing constraints at the design stage of the control algorithm. However, ensuring safety in an online fashion is challenging because the constraints may be unknown. Safe learning algorithms, such as SafeOpt from \cite{Safe_Sui2015}, use surrogate models to find a safe controller. However, safe learning algorithms based on surrogate models can become computationally expensive. In this paper, we extend the SafeOpt algorithm from \cite{Safe_Sui2015} to improve the performance of the algorithm and increase flexibility.

The SafeOpt algorithm proposed by \cite{Safe_Sui2015} is an iterative algorithm that uses Gaussian processes to learn unknown functional form of the objective and the constraints. SafeOpt ensures that in every iteration the new points that potentially maximize the objective while satisfying the constraints are chosen based on confidence intervals from the Gaussian processes. The choice of the new points is done by analysing the safety of a selected number of points from the whole search space. This approach is appealing because it does not rely on derivative information, which may often be unavailable. However, looking at the whole search space is equivalent to performing an exhaustive search and can quickly become computationally expensive, as indicated by \cite{BayesianBerkenkamp2021}. We reformulate the exhaustive search as a series of optimization problems to find the next recommended points. 

A review of safe learning methods related to SafeOpt was done by \cite{Safe_Kim2021}. Most of these algorithms also use discretized search space to find the optimum. However, already \cite{Safe_Fiducioso2019} identified the issue of computational effort related to looking at the whole search space in SafeOpt. The idea of merging derivative-free optimization with safe learning was explored by \cite{Constrained_Duivenvoorden2017} who proposed to find the recommended point by solving auxiliary optimization problems with particle swarm methods. They preserved the idea of SafeOpt to use confidence intervals of Gaussian processes in every iteration, but redefined the way of choosing new points to make it suitable for particle swarm methods. Due to the heuristic nature of particle swarm methods, their approach needed auxiliary adjustments to ensure good performance of the swarm. The reformulation we are proposing in the current paper allows avoiding heuristics while preserving the way the new points are chosen in every iteration of SafeOpt, at the same time improving the computational performance of the algorithm.

To preserve the character of SafeOpt, we solve the reformulated optimization problems using direct search methods. Similarly to particle swarm methods,  direct-search methods belong to derivative-free methods. In contrast to particle swarm methods, direct search provides explicit stopping criteria with convergence guarantees \cite[Ch. 2.]{Derivative_Audet2017}. In the current paper, we exploit the properties of pattern search to emulate the discretized structure of the search space required by SafeOpt. 

The contributions of the paper are twofold:
\begin{itemize}
    \item We reformulate SafeOpt as a series of optimization problems preserving the way new points are chosen in every iteration;
    \item We enable finding a trade-off between the computational time and accuracy of the solution by using direct search optimization methods.
\end{itemize}
We demonstrate with numerical examples that the proposed reformulation mitigates the computational effort while preserving safety.

The rest of the paper is structured  as follows. Section \ref{sec:Background} presents the necessary background and introduces the SafeOpt algorithm. The reformulation of the algorithm as a series of optimization problems is shown in Section \ref{sec:Reformulation}. Numerical examples in Section \ref{sec:Examples} show the performance of the reformulated algorithm. The paper ends with discussion in Section \ref{sec:Conclusions}.

\section{Background}
\label{sec:Background}
\subsection{Problem formulation}
In its basic form, SafeOpt from \cite{Safe_Sui2015} is designed to solve problems of the form:
\begin{subequations} \label{eqn:SafeOpt}
\begin{align}
\max_{x}& \quad f(x) \label{eq:ObjRewr}\\
\text{subject to}    &\quad g_j(x)\geq 0,~ j=1,\ldots,J\label{eq:CstrRewr},\\
    &\quad x\in\mathcal{A}. \label{eq:SetRewr}
\end{align}
\end{subequations}
where $x\in\mathcal{A}\subset\R^n$ is a vector of decision variables, $f:\R^n\rightarrow \R$ is the objective function to be maximized, $g_j:\R^n\rightarrow \R$ is one of $J$ constraints that must be satisfied. It is assumed that the functional form of neither $f$ nor $g_j$ is known, but there exist oracles that can provide noise-corrupted values of $f$ and $g_j$. As a result, they can be treated as outputs, whereas $x$ can be treated as an input in Gaussian process regression.

\subsection{Gaussian processes}
Following \cite{SafeBerkenkamp2016}, we use Gaussian processes to approximate both the objective and the constraints using measurements, i.e. we find approximations $J_i(x):\mathcal{A}\rightarrow \R$ where $i=0$ corresponds to the objective function, $i=1,\ldots,J$ corresponds to the constraints \eqref{eq:CstrRewr}. Gaussian process regression assumes that the values $J(x_0),J(x_1),\ldots,J(x_P)$ corresponding to different $x$ are random variables, with joint Gaussian distribution for any finite $P$. The prior information about the functions $J_i$ is defined by the mean and the covariance function $k(x_i,x_j)$.

Assume that we have access to noisy measurements $\hat{J}_i(x)=J(x)+\omega$ where $\omega\sim \mathcal{N}(0,\sigma^2_{\omega})$. To use Gaussian processes corresponding to $J_i$, we need to predict the value of $J_i$ at an arbitrary point $\hat{x}$ using only $R$ past measurement data $\mathbf{J}_i=[\hat{J}_i(x_r)]_{r=1,\ldots,R}$. From  \cite{SafeBerkenkamp2016}, the mean and the variance of the prediction are:
\begin{equation}
    \mu_i(\hat{x}) = \mathbf{k}_R(\hat{x})(\mathbf{K}_R+\mathbf{I}_R\sigma^2_{\omega})^{-1} \mathbf{J}_i
    \label{eq:mean}
\end{equation}
\begin{equation}
    \sigma^2_{R,i}(\hat{x}) = k(\hat{x},\hat{x})-\mathbf{k}_R(\hat{x})(\mathbf{K}_R+\mathbf{I}_R\sigma^2_{\omega})^{-1} \mathbf{k}_R^{\T}(\hat{x})
        \label{eq:variance}
\end{equation}
where $\mathbf{J}_i$ is a vector of $R$ observed noisy values, $i=0,\ldots,J$, the matrix $\mathbf{K}_R$ contains the covariances of past data, $k(x_a,x_b)$, $a,b=1,\ldots,R$, and $\mathbf{k}_R(\hat{x})$ contains the covariances between the new point and the past data. The identity matrix of size $R\times R$ is denoted by $\mathbf{I}_R$. 

The mean and the variance are then used to find the lower and upper confidence bounds:
\begin{equation}
\label{eq:Bounds}
\begin{aligned}
    l_R(x,i)=&{}\mu_i(x)-\beta\sigma_{R,i}(x),\\
    u_R(x,i)=&{}\mu_i(x)+\beta\sigma_{R,i}(x)
    \end{aligned}
\end{equation}
where $\beta$ corresponds to the desired confidence level.

\subsection{SafeOpt}

\subsubsection{Default SafeOpt}
\label{sec:DefaultSafeOpt}
The current paper is based on the version of SafeOpt proposed by \cite{SafeBerkenkamp2016}. The algorithm requires an \emph{initial safe set} set $S_0$ containing a number of initial points, such that the constraints \eqref{eq:CstrRewr} are satisfied. Using $S_0$, from \eqref{eq:Bounds}, we obtain the upper and lower bound corresponding to the initial safe set. The safe set in iteration $n$ is defined as: 
\begin{equation}
       S_n=\bigcap_{i=1,\ldots,J}\lbrace x\in A:l_n(x,i)\geq J_{\min} \rbrace
       \label{eq:SafeSetAlg}
\end{equation}
where the parameter $J_{\min}\geq 0$ is a design parameter that defines safety \citep{SafeBerkenkamp2016} and $A\subset\mathcal{A}$ is a discrete search space. The algorithm then looks for a recommended point $x_n\in S_n$ by finding a trade-off between \emph{maximizers} corresponding to the potential optimum and \emph{expanders} related to extending the current safe set $S_n$. 

The set of maximizers $M_n$ is defined as:
    \begin{equation}
        M_n=\lbrace x\in S_n:u_n(x,0)\geq \max_{x\in S_n}l_n(x,0) \rbrace.
        \label{eq:Maximizers}
    \end{equation}

\cite{SafeBerkenkamp2016} define the set of expanders $G_n$ as:
\begin{equation}
\label{eq:Expanders}
    G_n=\lbrace x\in S_n: |\mathcal{G}(x)|>0  \rbrace
\end{equation}
where $|\cdot|$ represents the cardinality of a set and 
\begin{equation}
\label{eq:AuxFcn}
    \mathcal{G}(\overline{x})=\lbrace x'\in A\setminus S_n:\forall j\; l_{n,j,(\overline{x},u_n(\overline{x},j))}(x')\geq J_{\min}\rbrace.
\end{equation}
where $l_{n,(\overline{x},u_n(\overline{x}))}(x)$ the lower bound for the point $x$ obtained from the auxiliary GP calculated for a given point $\overline{x}\in S_n$. The auxiliary GP is created assuming that we observed the optimistic upper bound $u_n(\overline{x},i)$ \citep{SafeBerkenkamp2016}. 
    
The new point is chosen as:
    \begin{equation}
        x_n=\underset{x\in G_n\cup M_n}{\text{argmax}}\max_{i} w_n(x,i)
        \label{eq:RecommendedPoint}
    \end{equation}
where $w_n(x,i)=u_n(x,i)-l_n(x,i)$. The point $x_n$ is then applied to the system to collect new observations about the constraints and the objective, and update the Gaussian processes using \eqref{eq:mean} and \eqref{eq:variance}.

The optimum after $n$ iterations is found as:
\begin{equation}
    x^*=\argmax_{x\in S_n}l_n(x,0).
    \label{eq:MaxEstim}
\end{equation}
Analysis of convergence and safety guarantees of the algorithm was done by \cite{SafeBerkenkamp2016}.

\subsubsection{Discretization of the search space}
Equation \eqref{eq:RecommendedPoint} describes searching for a recommended point $x_n$ in the search space defined by the union of the set of expanders $G_n$ and maximizers $M_n$. The search space $G_n\cup M_n$ is always bounded if the overall search space $\mathcal{A}$ is bounded. If the search space $G_n\cup M_n$ is also discrete, the value of $x_n$ from \eqref{eq:RecommendedPoint} can be found by exhaustive search in the set $G_n\cup M_n$. The number of discretization points $N$ will then define how big the search space for the exhaustive search is. Because of the need for an exhaustive search, finding $x_n$ becomes a challenge for large values of $N$ because it imposes a significant computational burden \citep{BayesianBerkenkamp2021}. 

The computational effort is primarily spent on finding the auxiliary Gaussian process to obtain $l_{n,j,(\overline{x},u_n(\overline{x},j))}(x')$ for all $x'\in A\setminus S_n$. Therefore, limiting the number of points $x'$ has the potential of mitigating the computational effort. One way of limiting the number of points relies on choosing a small $N$ corresponding to a coarse discretization of $\mathcal{A}$ so that the number of points in the set $A\setminus S_n$ is small. However, a small $N$ corresponds also to the low accuracy of the solution to \eqref{eqn:SafeOpt}. To avoid the choice of discretization and get more flexibility in setting the accuracy, we reformulate the SafeOpt algorithm for finding the recommended value from \eqref{eq:RecommendedPoint} as a series of optimization problems while preserving the definitions of the sets of maximisers \eqref{eq:Maximizers} and the expanders \eqref{eq:Expanders}.

\subsection{Pattern search algorithm}
Pattern search methods belong to the group of direct search optimization methods and rely on evaluating a number of candidate points around a selected point, which are chosen following a given \emph{pattern}.

Given a selected point $x^k\in\R^n$, a \emph{pattern} defines a set of vectors in $\R^n$ where the algorithm looks for an incumbent point. Formally, the algorithm uses a \emph{mesh}, defined as \citep{Derivative_Audet2017}:
\begin{equation}
    \mathcal{M}^k:=\lbrace x^k+\delta^k Dy:y\in \N^p \rbrace
    \label{eq:Mesh}
\end{equation}
where $\delta^k>0$ describes the mesh size of $\mathcal{M}^k$, and $D=GZ$ where $G\in\R^{n\times n}$ is an invertible matrix, $Z\in \R^{n\times p}$ is such that the columns of $Z$ form a positive spanning set in $\R^n$. Let us denote with $\mathbb{D}$ the columns of $D$. A positive spanning set $\mathbb{D}^k\subseteq \mathbb{D}$ is called a \emph{pattern}. An in-depth description of the algorithm was provided by \cite{Derivative_Audet2017}.

\subsubsection{Mesh and discretization}
Pattern search evaluates the points defined by the current mesh to find $x^{k+1}$ that improves the value of the objective function while satisfying constraints \citep{Globally_Lewis2002}. Then the mesh size is adjusted, $\delta^{k+1}=2\delta^k$. However, if it is impossible to find a point $x^{k+1}$ for the current mesh size, the algorithm decreases the mesh size $\delta^{k}=0.5\delta^k$ and evaluates the points defined by the mesh with the updated size. 

\begin{figure}
\psfrag{DISCRETISED}[][]{\scriptsize{\textsf{Discretized}}}
\psfrag{SPACE}[][]{\scriptsize{\textsf{\hspace{0.3cm}space $A$}}}
\psfrag{SEARCH SPACE}[][]{\scriptsize{\textsf{\hspace{-0.3cm}Search space $\mathcal{A}$}}}

\psfrag{MESH}[][]{\scriptsize{\textsf{Meshes}}}
\psfrag{MA}[][]{\tiny{\textsf{$\mathcal{M}_A$}}}
\psfrag{MB}[][]{\tiny{\textsf{$\mathcal{M}_B$}}}

     \centering
         \includegraphics[width=0.45\textwidth]{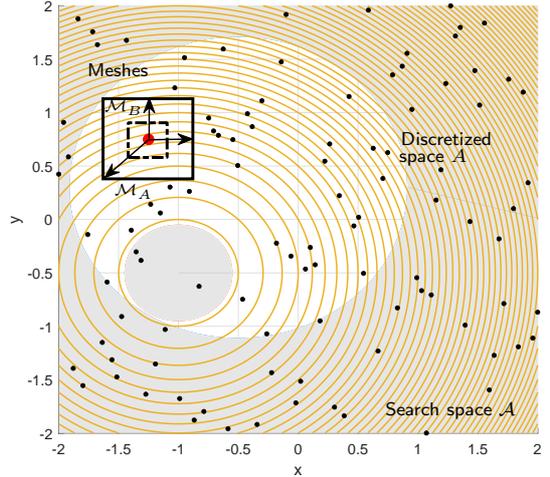}
        \caption{Illustration of the discretized search space (black dots) and the meshes used by pattern search with the pattern $[0,1], [1,0], [-1,-1]$ (arrows) around the current point (red circle). Mesh $\mathcal{M}_B$ has been obtained by decreasing the size of mesh $\mathcal{M}_A$ by 0.5}
        \label{fig:MeshDiscrete}
\end{figure}

The mesh defined by \eqref{eq:Mesh} can be intuitively understood as local discretization around the current point $x^k$. An illustration is provided for two-dimensional search space in Fig. \ref{fig:MeshDiscrete}. The default SafeOpt uses the discretized search space $A\subset\mathcal{A}$ to evaluate safety (black dots). The squares show the meshes used by pattern search with the pattern $[0,1], [1,0], [-1,-1]$ around the current point (red circle). The mesh $\mathcal{M_B}$ (dashed) has been obtained by decreasing the size of mesh $\mathcal{M}_A$ (solid) by 0.5.

The algorithm stops if the mesh size becomes smaller than a given threshold $\delta^k\leq\varepsilon$. We explore the flexibility provided by using the mesh size as a stopping criterion in pattern search to introduce new stopping criteria for the reformulated SafeOpt algorithm. We use the norm of the solution as a stopping criterion to achieve the desired accuracy. The impact of the new criterion will be shown in Section \ref{sec:Examples}. 

\section{SafeOpt as optimization}
\label{sec:Reformulation}
We now present the main result of this paper. We solve two separate optimization problems to find the new recommended value $x_n$:
\begin{equation}
    P_1: \quad \underset{x\in M_n}{\text{max}} \max_{i} w_n(x,i),
        \label{eq:ProblemCast2}
\end{equation}
\begin{equation}
    P_2: \quad \underset{x\in G_n}{\text{max}} \max_{i} w_n(x,i)
    \label{eq:ProblemCast1}
\end{equation}
The new recommended value $x_n$ is obtained as the point from $G_n\cup M_n$ such that $w(x_n)=w_{\max}$ with
    \begin{equation}
        w_{\max}:={\text{max}}\Big\lbrace\underset{x\in G_n}{\text{max}} w(x), \underset{x\in M_n}{\text{max}} w(x)\Big\rbrace.
        \label{eq:MaxW}
    \end{equation}
where $w(x):=\max_{i} w_n(x,i)$.  To solve the two problems $P_1$ and $P_2$ using numerical solvers, we need now to rewrite the search space of each problem in the form of constraints. 

\subsection{Reformulation of the maximizers}
\label{sec:Maximizers}
From the definition of the safe set from \eqref{eq:SafeSetAlg}, we obtain that:
\begin{equation}
    x\in S_n \iff x\in\mathcal{A} \text{ and }  \forall j=1,\ldots,J\quad l_n(x,j)\geq J_{\min}.
    \label{eq:SafeSetIneq}
\end{equation}
From the definition of the maximisers \eqref{eq:Maximizers} we obtain:
\begin{equation}
    x\in M_n \iff x\in S_n \text{ and } u_n(x,0)\geq l^*
    \label{eq:MaxDef}
\end{equation}
where 
\begin{subequations} \label{eq:InternalSafe}
\begin{align}
l^* =  \max_{z}&\;l_n(z,0)\\
\text{subject to } &  l_n(z,j)\geq J_{\min},~ \forall j =1 ,\ldots,J.
\end{align} 
\end{subequations}

We note in \eqref{eq:ProblemCast2} that $w_n(\cdot,i)$, $w_n(\cdot,j)$ are independent from each other for $i\neq j$. Therefore, the objective function \eqref{eq:ProblemCast2} can be reformulated yielding $J$ separate problems $P_1^k$, $k=1,2,\ldots,J$. Using \eqref{eq:SafeSetIneq} and \eqref{eq:MaxDef}, we obtain:
\begin{subequations} \label{eq:ReformulatedMaxAll}
\begin{align}
P_1^k:\quad \underset{y\in \mathcal{A}}{\text{max}}&\quad w_n(x,k)\\
\text{subject to}    &\quad l_n(x,j)\geq J_{\min},~  \forall j=1,\ldots,J,\\
&\quad  u_n(x,0)\geq l^*,
\end{align}
\end{subequations}
where $l^*$ is obtained from \eqref{eq:InternalSafe}. The problem in  \eqref{eq:InternalSafe} is independent of $x$ and can be solved separately.

Let us denote a solution to the problem $P_1^k$ as $x_1^{k*}$. Then the solution to \eqref{eq:ProblemCast1} is found as:
\begin{equation}
    x_1^*=\max_{k=1,\ldots,J} x_1^{k*},
\end{equation}
obtained for $k_1^*$.

\subsection{Reformulation of the expanders}
\label{sec:Expanders}
In this section, we describe the proposed reformulation of \eqref{eq:ProblemCast1}. From \eqref{eq:Expanders} we get:
\begin{subequations}
\label{eq:ExpandersRef}
\begin{align}
\underset{x}{\text{max}}&\quad  \max_{i} w_n(x,i)\\
\text{subject to}    &\quad x\in S_n,\\
&\quad |\mathcal{G}(x)|> 0
\end{align}
\end{subequations}
where $\mathcal{G}(\cdot)$ is given by \eqref{eq:AuxFcn}. From \eqref{eq:AuxFcn} we notice that $|\mathcal{G}(x)|> 0$ \ if there exists at least one point $x'\in \mathcal{A}\setminus S_n$ such that the condition:
\begin{equation}
    \forall j\; l_{n,j,(x,u_n(x,j))}(x')\geq J_{\min}
\end{equation}
is satisfied. Thus, we obtain:
\begin{subequations}
\label{eq:ExpandersRefInf}
\begin{align}
\underset{x,x'}{\text{max}}&~ \max_i w_n(x,i)\\
\text{subject to}    &~ x\in S_n, \\
&~  l_{n,j,(x,u_n(x,j))}(x')\geq J_{\min},~  \forall j=1,\ldots,J,\\
&~  x'\in \mathcal{A}\setminus S_n.
\end{align}
\end{subequations}
From the definition of the safe set \eqref{eq:SafeSetAlg}, we get that:
\begin{equation}
    x'\in \mathcal{A}\setminus S_n \iff  x'\in \mathcal{A} \text{ and } \exists k: l_n(x',k)< J_{\min}
\end{equation}
Then we have the following equivalence:
\begin{equation}
  \exists k: l_n(x',k)< J_{\min} \iff \min_s l_n(x',s)<J_{\min}  
  \label{eq:SafeEquivalence}
\end{equation}
Then we obtain:
\begin{subequations}
\label{eq:ExpandersRefInfOpt}
\begin{align}
\underset{x,x'}{\text{max}}&\quad \max_{i} w_n(x,i)\\
\text{subject to}    &\quad  \min_s l_n(y,s) \geq J_{\min}, \\
&\quad  \min_{s}l_{n,s,(x,u_n(x,s))}(x')\geq J_{\min}\label{eq:ExpandersSetIneq},\\
&\quad  \min_s l_n(x',s) < J_{\min}.
\end{align}
\end{subequations}
If there exists a solution to \eqref{eq:ExpandersRefInfOpt}, the optimization problem can be solved using a derivative-free method, such as pattern search. However, the problem \eqref{eq:ExpandersRefInfOpt} may be infeasible if the set of expanders $G_n$ is empty. To avoid potential infeasibility, we relax \eqref{eq:ExpandersRefInfOpt} as:
\begin{subequations}
\label{eq:ExpandersRefFeas}
\begin{align}
\underset{x,x'}{\text{max}}&\quad q(x,x')\\
\text{subject to}    &\quad  \min_s l_n(x,s) \geq J_{\min}, \label{eq:Safey}\\
&\quad  \min_s l_n(x',s) < J_{\min} 
\end{align}
\end{subequations}
where 
\[
\begin{aligned}
q(x,x')=&
\max_{i} w_n(x,i) - \\
&\sigma\min\lbrace 0,\min_s l_{n,s,(x,u_n(x,s))}(x')-J_{\min} \rbrace
\end{aligned}
\]
where $\sigma>0$ enables trading off feasibility and optimality. The problem from \eqref{eq:ExpandersRefFeas} is feasible if $S_n\subset \mathcal{A}$ in the strict sense, i.e. $S_n\neq \mathcal{A}$. Detecting infeasibility is out of the scope of most derivative-free solvers so we ensured that the problem is feasible in the relaxation \eqref{eq:ExpandersRefFeas}.

Doing the same reformulation as in \eqref{eq:ReformulatedMaxAll}, we get:
\begin{subequations} \label{eq:ReformulatedExp}
\begin{align}
P_2^k:\quad \underset{x,x'\in \mathcal{A}}{\text{max}}&\quad q_k(x,x')\\
\text{subject to}    &\quad  \min_s l_n(x,s) \geq J_{\min},\\
&\quad  \min_s l_n(x',s) < J_{\min} 
\end{align}
\end{subequations}
where 
\begin{align}
q_k(x,x')&=w_n(x,k) - \notag \\
    &\min\lbrace 0,\min_s \lbrace l_{n,s,(x,u_n(x,s))}(x')-J_{\min} \rbrace\rbrace.
\label{eq:ExpRefFeas}
\end{align}
Let us denote a solution to the problem $P_2^k$ as $x_2^{k*}$. Then the solution to \eqref{eq:ProblemCast2} is found as:
\begin{equation}
    x_2^*=\argmax_{k=1,\ldots,J} w_n(x_2^{k*},k).
\end{equation}
obtained for $k_2^*$. Then from \eqref{eq:MaxW} we get:
\begin{equation}
    x_n =\argmax_{\{x_1^*,x_2^*\}}\lbrace w_n(x_2^*,k_2^*),w_n(x_1^*,k_1^*) \rbrace.
    \label{eq:Maxw2}
\end{equation}
The reformulation proposed in Sections \ref{sec:Maximizers} and \ref{sec:Expanders} is summarized in Algorithm \ref{alg:SafeOptRef}. 

We note that the problems \eqref{eq:ReformulatedMaxAll}, \eqref{eq:ReformulatedExp} use the same definitions of the maximisers and the expanders as \cite{SafeBerkenkamp2016}. In particular, the proposed reformulation is independent from the chosen pattern search solver. Thus, the problems \eqref{eq:MaxEstim}, \eqref{eq:ReformulatedMaxAll}, \eqref{eq:ReformulatedExp} (lines 5, 6, and 7 in Algorithm \ref{alg:SafeOptRef}) are general and can be solved with other methods. For instance, if the problems are differentiable, the reformulation enables using gradient-based methods.

\setlength{\algomargin}{1.1em}
\begin{algorithm2e}[t]
\SetAlgoLined
  \KwInput{Initial safe set $S_0=\lbrace x^0,x^1,\ldots,x^K\rbrace\subset\mathcal{A}$, desired tolerances $\epsilon_1$, $\epsilon_2$, maximal number of iterations $M$, desired confidence limit $\alpha$, desired safety threshold $J_{\min}$}
  \KwOutput{Optimal solution $x^*$}
  Set $n\leftarrow 1$, compute $F_n=\lbrace f(x^i) \rbrace_{i=1,\ldots,K}$, $G_{k,j}=\lbrace g_j(x^i) \rbrace_{i=1,\ldots,K}$ for $j=1,\ldots,J$, set $S_n\leftarrow S_0$.
  
  \Repeat{$n\!\geq\! M\!\And\!\|x_n^r\!-\!x_{n-1}^r\|\!\leq\! \epsilon_1 \!\And\! \|f(x_n^r)\!-\!f(x_{n-1}^r)\|\!\leq\! \epsilon_2$}{
Using $S_n$ and $F_n$, find a Gaussian process $GP_f$ with lower bounds $l_n(x,0)$, upper bounds $u_n(x,0)$, 

Using $S_n$ and $G_{k,j}$ find $J$ Gaussian processes $GP_{g,j}$ with lower bounds $l_n(x,j)$, upper bounds $u_n(x,j)$

Solve \eqref{eq:MaxEstim}, obtaining $x^*_n$ and $l^*=l_n(x^*_n,0)$

Solve $P^j_1$ for all $j=1,\ldots,J$ from \eqref{eq:ReformulatedMaxAll}

Solve $P^j_2$ for all $j=1,\ldots,J$ from \eqref{eq:ReformulatedExp}

\eIf{$q_j(x_2^{j*},x^{'*})\leq w_n(x_2^{j*},j)$}{
Set $x_n^r\leftarrow x_1^*$}{  Solve \eqref{eq:Maxw2} and set $x_n^r\leftarrow \argmax_{\{x_1^*,x_2^*\}}\lbrace w_n(x_2^*,k_2^*),w_n(x_1^*,k_1^*) \rbrace$

}

Set $n\leftarrow n+1$, set $F_n\leftarrow F_{n-1}\cup\lbrace f(x^r) \rbrace$, $G_{k,j}\leftarrow G_{n-1}\cup\lbrace g_j(x^r) \rbrace$ for $j=1,\ldots,J$, set $S_n\leftarrow S_{n-1}\cup \lbrace x^r \rbrace$

}
 \caption{Reformulated SafeOpt\label{alg:SafeOptRef}}
\end{algorithm2e}

\section{Examples}
\label{sec:Examples}
This section presents the performance of the proposed reformulation. All tests were performed in Windows 10, using Matlab 2021a on a laptop with an AMD Ryzen 7 PRO 5850U, 8 cores, with 32\,GB of RAM.

\subsection{Safety with non-convex feasible set}
\begin{figure}
     \centering
         \includegraphics[width=0.45\textwidth]{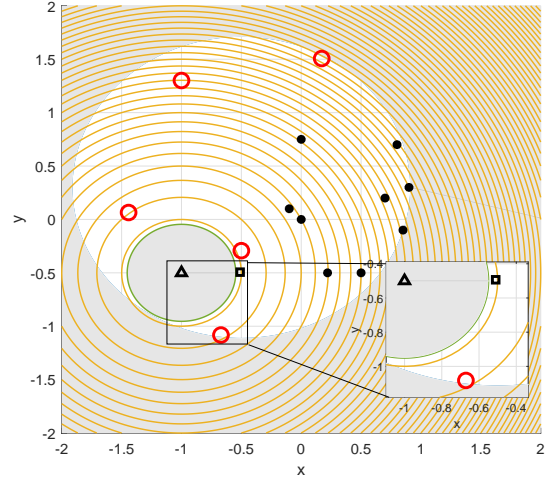}
        \caption{Illustration of \eqref{eqn:QuadProblem}. The feasible set is inside the white area. The unconstrained maximum of the quadratic function \eqref{eq:Quad} is denoted with a triangle and the found maximum is indicated with a square. The algorithm used the recommended points shown with solid circles obtained starting from the safe set denoted with black dots }
        \label{fig:QuadProblem}
\end{figure}

We show the safety of the proposed reformulation for solving an optimization problem with a non-convex feasible set:
\begin{subequations} \label{eqn:QuadProblem}
\begin{align}
\max_{x}& \quad -(x+1)^2-(y+0.5)^2 \label{eq:Quad}\\
\text{subject to}    &\quad 2-(x_1+0.5)^2-(y-0.3)^2\geq 0,\\
&\quad (x+1)^2+(y+0.5)^2-0.2\geq 0. \label{eq:CstrOut}
\end{align}
\end{subequations}
The maximum of \eqref{eq:Quad} is obtained for $x=-1$, $y=-0.5$, which lies outside the feasible set. A feasible maximum is anywhere on the boundary of the set defined by \eqref{eq:CstrOut}. We set $\epsilon_1=\epsilon_2=0.001$ and $J_{\min}=0$. An illustration is shown in Fig. \ref{fig:QuadProblem}. The thin circles denote level sets of the cost \eqref{eq:Quad}. The unconstrained maximum is denoted with a triangle. The constrained optima are on the level set shown with a green line. The point found by the algorithm was at $x=-0.51$ and $y=-0.5$ (square in Fig. \ref{fig:QuadProblem}). The red circles indicated the interim recommended points if the algorithm was initialized with the safe set shown with the black dots. All of the interim recommended points are inside the feasible set which confirms that safety was preserved. In particular, the inset in the bottom right corner shows a recommended point close to the boundary of the feasible set.

\subsection{Unknown constraints}
In this example, we apply the proposed reformulation to a problem of tuning a cascade PID controller for ball-screw drive from \cite{Cascade_Khosravi2020,Performance_Khosravi2022}. 

\subsubsection{Problem setup}
\begin{figure}
\psfrag{speed}[][]{\scriptsize{$S$}}
\psfrag{position}[][]{\scriptsize{$P$}}
\psfrag{speedsp}[][]{\scriptsize{$S_s$}}
\psfrag{positionsp}[][]{\scriptsize{$P_s$}}
\psfrag{speedctr}[][]{\scriptsize{$C_s(s)$}}
\psfrag{positictr}[][]{\scriptsize{$C_p(s)$}}
\psfrag{system}[][]{\scriptsize{$G(s)$}}
\psfrag{transfer}[][]{\scriptsize{$1/s$}}
     \centering
         \includegraphics[width=0.45\textwidth]{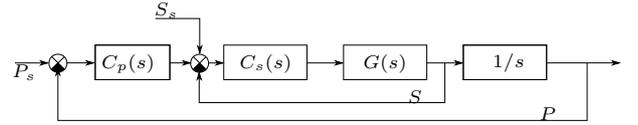}
        \caption{Block diagram of a ball–screw drive with transfer function $G(s)$. The objective is to follow the position set point $P_s$ ensured by a proportional controller $C_p(s)$ in cascade with a speed controller $C_s(s)$. }
        \label{fig:BlockDiagram}
\end{figure}

The control structure is shown in Fig. \ref{fig:BlockDiagram} and the parameters are the same as used by \cite{Performance_Khosravi2022}. We want to find a parameter $K_p$ for the position controller $C_p(s)$ and the parameters $K_v$ and $K_{vi}$ for the speed control $C_s(s)$ to minimize the integrated absolute error in the position and the overshoot in the speed:
\begin{equation}
    J := \gamma\int\limits_0^{t_f} |P(\tau)-P_s(\tau)| \d \tau + \max_{\tau\in[0,t_f]} S(\tau)
\end{equation}
with $\gamma=1000$ putting emphasis on the position tracking the desired reference $P_s$.

\begin{figure*}
     \centering
     \begin{subfigure}[b]{0.44\textwidth}
         \centering
         \includegraphics[width=\textwidth]{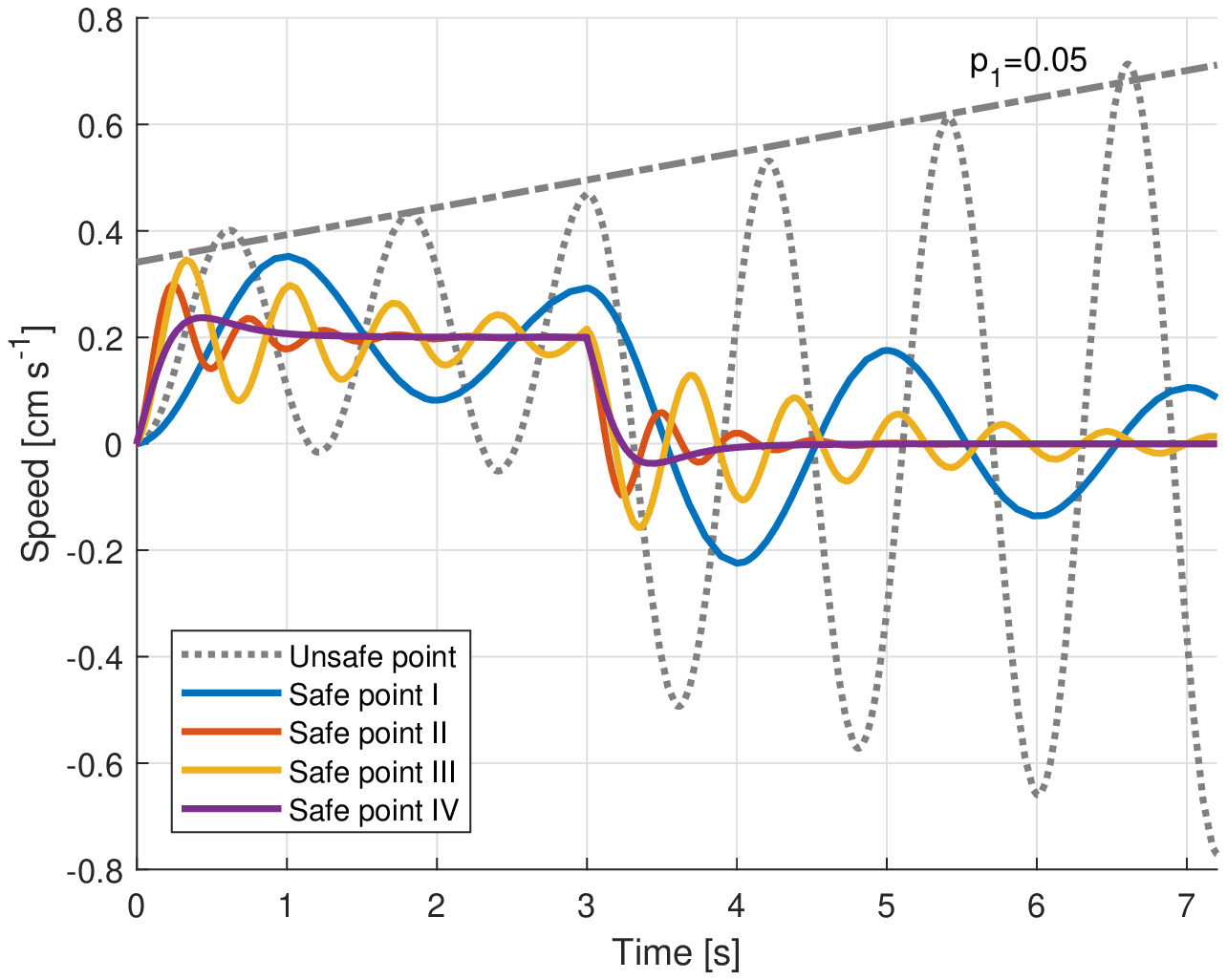}
         \caption{Speed for the initial safe set and the tuning set point}
         \label{fig:InitialSafe1}
     \end{subfigure}
     ~
     \begin{subfigure}[b]{0.44\textwidth}
         \centering
         \includegraphics[width=\textwidth]{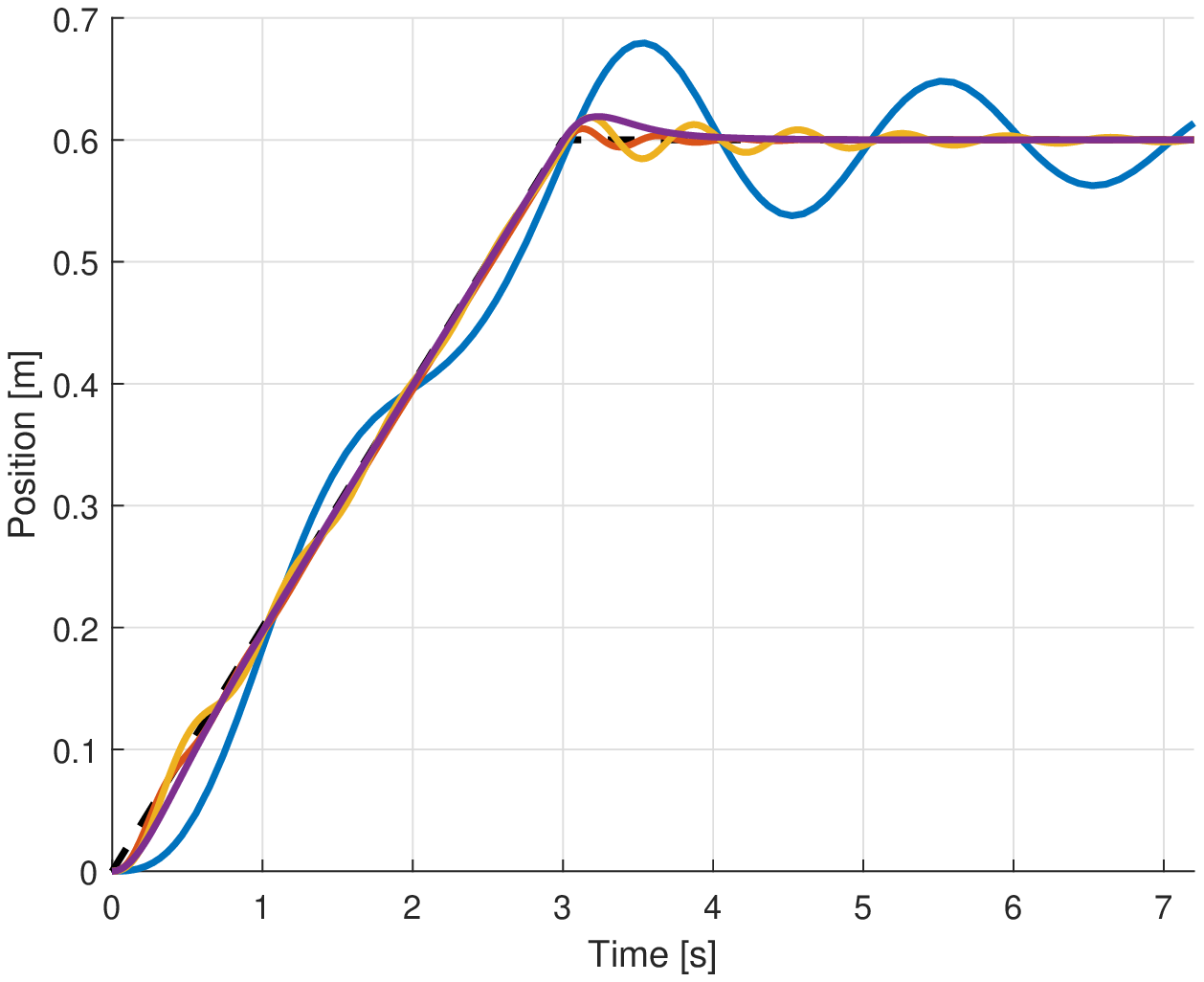}
         \caption{Position for the initial safe set and the tuning set point}
         \label{fig:InitialSafePos1}
     \end{subfigure}
     \hfill
          \begin{subfigure}[b]{0.44\textwidth}
         \centering
         \includegraphics[width=\textwidth]{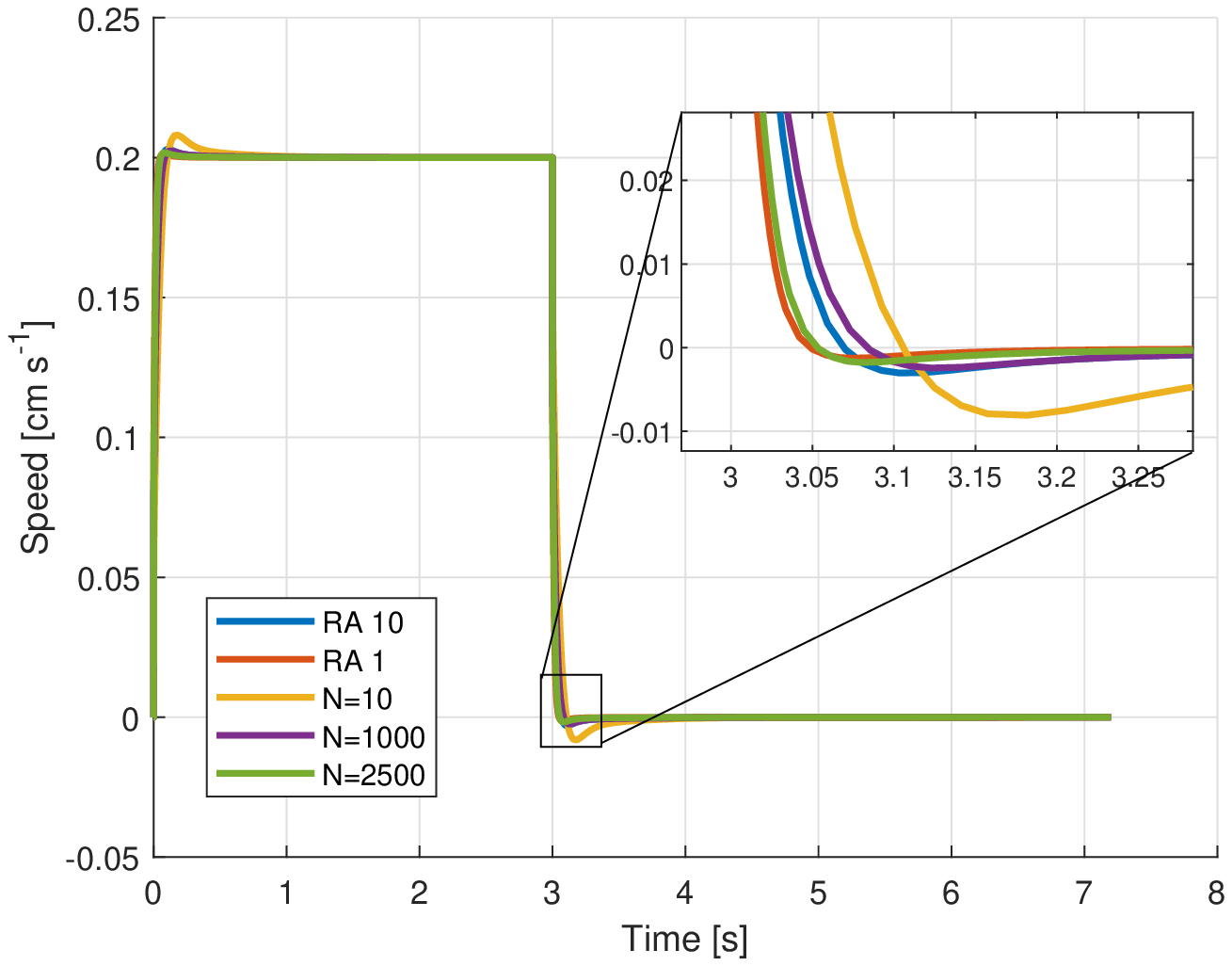}
         \caption{Speed for the solutions and the tuning set point}
         \label{fig:SpeedRes1}
     \end{subfigure}
     ~
     \begin{subfigure}[b]{0.44\textwidth}
         \centering
         \includegraphics[width=\textwidth]{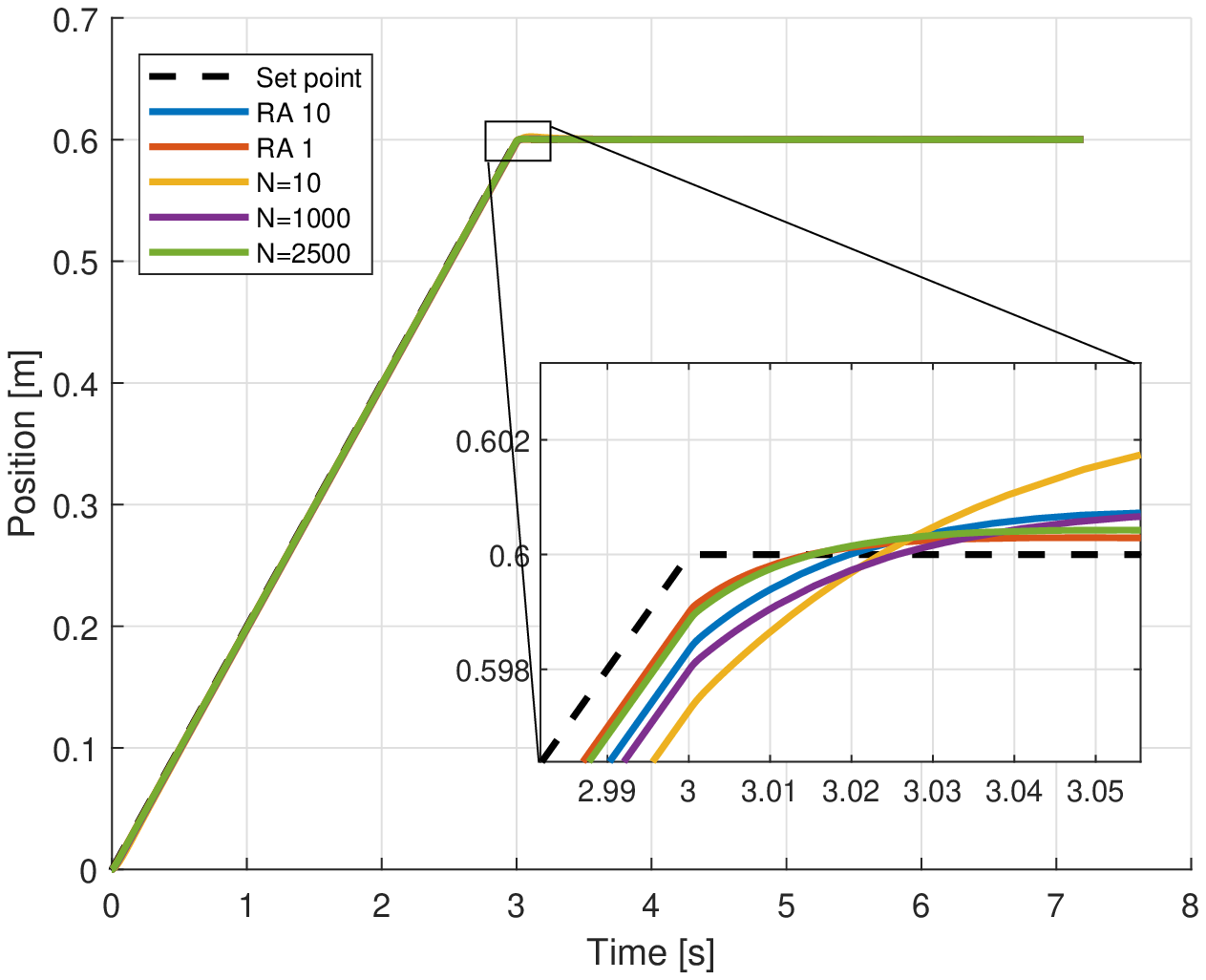}
         \caption{Position for the solutions and the tuning set point}
         \label{fig:PositionRes1}
     \end{subfigure}
          \hfill
          \begin{subfigure}[b]{0.44\textwidth}
         \centering
         \includegraphics[width=\textwidth]{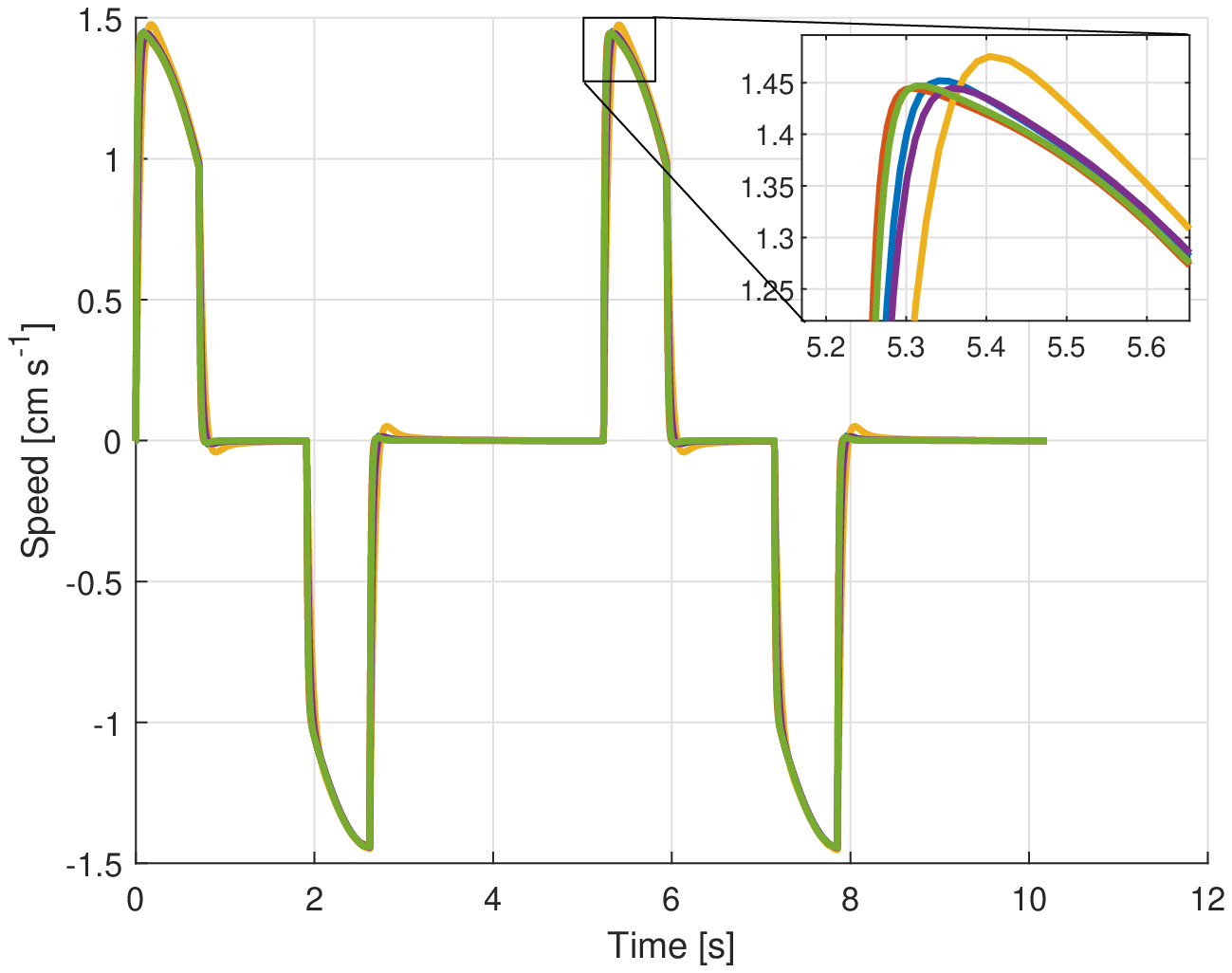}
         \caption{Speed for the solutions and the truncated sinusoidal set point}
         \label{fig:SpeedResSin1}
     \end{subfigure}
     ~
     \begin{subfigure}[b]{0.44\textwidth}
         \centering
         \includegraphics[width=\textwidth]{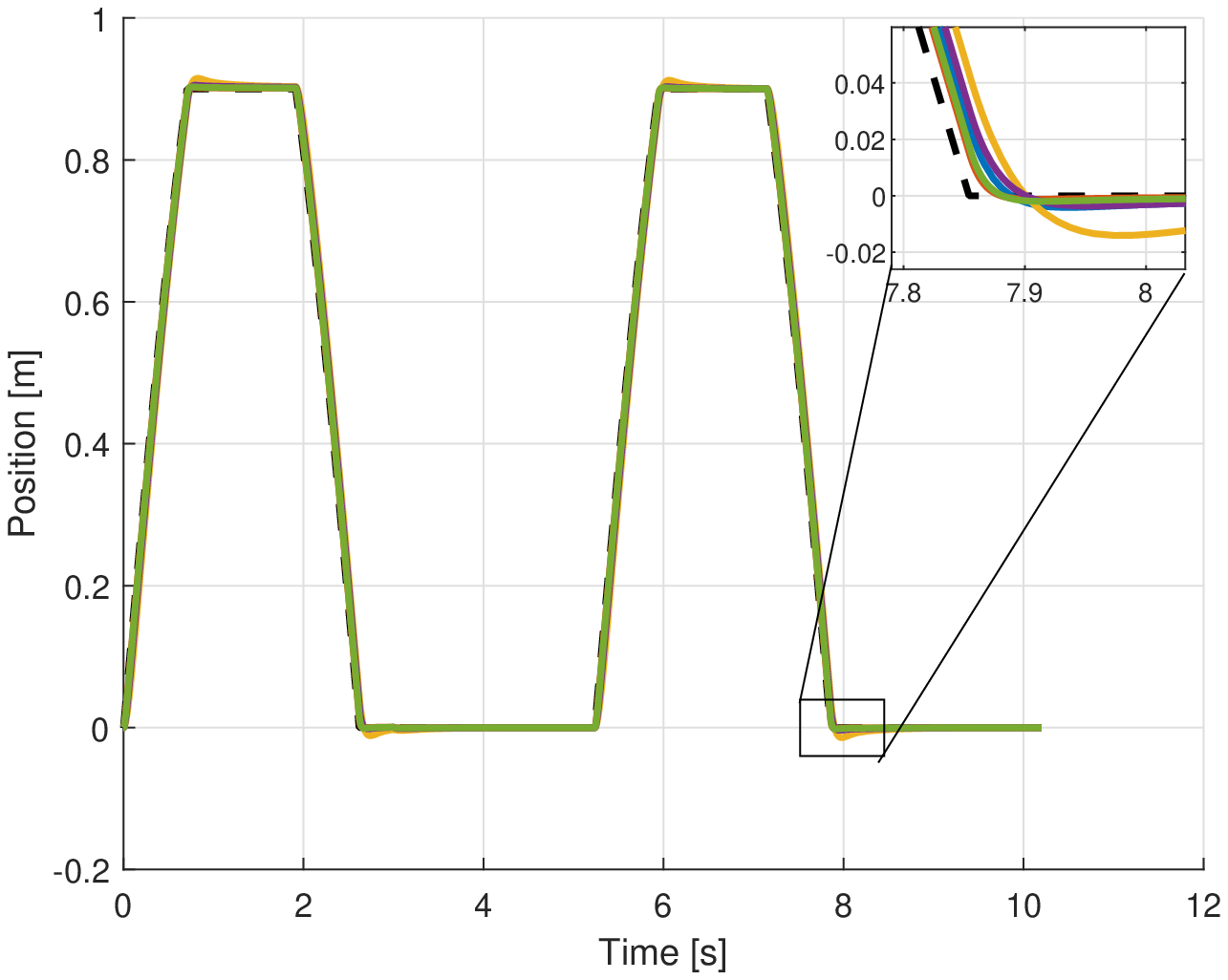}
         \caption{Position for the solutions and the truncated sinusoidal set point}
         \label{fig:PositionResSin1}
     \end{subfigure}
        \caption{Speed and position corresponding to the initial safe set $S_0$ with four points and to the resulting parameters applied to the training trajectory, and to the sinusoidal trajectory truncated at the position of 0.9 cm and zero}
        \label{fig:InitialSafeSet}
        
\end{figure*}

The system must satisfy a stability constraint. To emulate human-driven PID tuning based on visual assessment of responses of the system, we measure stability as the slope $p_1$ of the peaks of the response of the system, with positive values indicating instability \citep[Ch. 4.4]{Advanced_Astroem2006}. The constraint was formulated as:
\begin{equation}
    h(K_p,K_v,K_{vi}):=p_1-\sigma\leq 0
    \label{eq:fitp1}
\end{equation}
where $\sigma=0.005$ was chosen to ensure that a system with no peaks, i.e. $p_1=0$, yields a value inside the feasible set. The grey lines in Fig. \ref{fig:InitialSafe1} show the speed trajectory for an unstable point (dotted line) and the corresponding value $p_1$ used to evaluate \eqref{eq:fitp1}. The dash-dotted line with slope $p_1$ is the linear fit to the peaks of the unstable trajectory.

The search space $\mathcal{A}=[0,110]\times [0,50]^2$ was chosen so that it contains unstable values. The initial safe set contains four points collected in Table \ref{tbl:InitialSafeSet}. The initial values were found in simulation. The corresponding trajectories for speed and position are shown in Fig. \ref{fig:InitialSafeSet}. 

\begin{table}[!tbp]
\caption{Four points from the initial safe set $S_0$, and an unsafe point, together with the corresponding value of the objective function and the constraint}
\label{tbl:InitialSafeSet}
\begin{tabular}{@{}lccccc@{}}
\toprule
             & $K_p$ & $K_v$ & $K_{vi}$ & Obj. \eqref{eq:RefObjPID} & Const. \eqref{eq:RefCstrPID} ($p_1$) \\ \cmidrule{1-6}
Unsafe point   & 30    & 0     & 5        & -388      & -4.6 (0.05)         \\
Safe point I   & 10    & 0     & 5        & -241      & 4.8 (-0.04)       \\
Safe point II  & 20    & 0.4  & 50       & -20       & 5.4 (-0.05)       \\
Safe point III & 42    & 0.3   & 12       & -39       & 5.8 (-0.05)       \\
Safe point IV  & 90    & 0.5   & 1        & -26       & 0.005 (0)         \\ \bottomrule
\end{tabular}
\end{table}

Setting $x:=[K_p,K_v,K_{vi}]^{\T}$, we obtain the problem structure from \eqref{eqn:SafeOpt}:
\begin{subequations} \label{eqn:SafeOptPID}
\begin{align}
\max_{x}& \quad -J(x) \label{eq:RefObjPID} \\
\text{subject to}    &\quad -h(x)\geq 0 \label{eq:RefCstrPID},\\
    &\quad x\in\mathcal{A}.
\end{align}
\end{subequations}
Table \ref{tbl:InitialSafeSet} shows the values of the objective \eqref{eq:RefObjPID} and the constraint \eqref{eq:RefCstrPID} for the initial safe set $S_0$ and $J_{\min}=0$.

\subsubsection{Performance}
To show the computational flexibility of the proposed reformulation, we explore the properties of the pattern search algorithm. In particular, we show that by adjusting the selected stopping criteria of the pattern search algorithm we reduce computational complexity, measured by the time needed to find a solution while preserving good solution performance, measured by the value of the objective function.

\begin{table*}[!tbp]
\centering
\caption{Performance assessment based on stopping criteria enabled by the reformulation and stopping criteria of pattern search: mesh tolerance $\varepsilon$ (scaled by $10^{-4}$ in the table), initial mesh size $\delta^0$ with the best results marked in bold}
\label{tbl:Performance}
\begin{tabular}{@{}ccccccccccP{2cm}P{2cm}@{}}
\toprule
$\epsilon_1$ & $\epsilon_2$ & $\varepsilon$ [$\times 10^{-4}$] & $\delta^0$ & $K_p^*$ & $K_v^*$ & $K_{vi}^*$ & Obj. \eqref{eq:RefObjPID} & Iterations & Time [s]       &  \parbox[c]{\hsize}{\centering Time [s]  for solving \eqref{eq:ReformulatedMaxAll}} &  \parbox[c]{\hsize}{\centering Time [s]  for solving \eqref{eq:ReformulatedExp}} \\ \midrule
0.001         & 0.001        & 0.01      & 1          & 107.2   & 45.3    & 41.2       & -3.4                      & 25         & 268            & 0.1                                      & 4.2                                      \\
0.01         & 0.01         & 0.01      & 1          & 55      & 31.2    & 49.3       & -6.4                      & 5          & 44             & 0.4                                      & 6.4                                      \\
0.1         & 0.1        & 0.01      & 1          & 55      & 31.2    & 49.3       & -6.4                      & 4          & 35             & 0.6                                      & 4.5                                      \\
0.1         &0.1        & 0.01     & 5          & 73.1    & 41.1    & 49.6       & -4.8                      & 8          & 68             & 0.4                                      & 5.6                                      \\
0.1         & 0.1         & 0.01     & 10         & 60      & 25.4    & 50         & -5.8                      & 2          & 23             & 0.6                                      & 9.8                                      \\
0.1         & 0.1         & 0.01      & 20         & 42      & 40      & 50         & -8.2                      & 3          & 37             & 0.7                                      & 9.1                                      \\
0.1         & 0.1        & 1      & 10         & 60      & 25.4    & 50         & -5.8                      & 2          & 13             & 0.3                                      & 4.3                                      \\
0.1        & 0.1         & 100      & 10         & 60      & 25.4    & 50         & -5.8                      & 3          & $\mathbf{6}$ & 0.02                                     & 0.1                                      \\
0.1        & 0.1         & 100       & 1          & 100     & 43.4    & 41         & -$\mathbf{3.4}$         & 8          & 17             & 0.01                                     & 0.1                                      \\ \bottomrule
\end{tabular}
\end{table*}

The results are shown in Table~\ref{tbl:Performance}. The first row shows that strict stopping criteria allowed for reaching a low value of the objective function. However, the overall time needed to get a solution is significant. This is because the algorithm keeps running until the solutions from two consecutive iterations are close. Relaxation of the stopping criteria enabled by the reformulation impacts primarily the number of iterations and thus the overall time to get a solution, reducing it from 268 s (4 min 45 s) to 35 s. Stopping the algorithm after a smaller number of iterations led to a larger value of the objective function, which increased from 3.4 to 6.4. Therefore, we see that in the proposed reformulation we can use the stopping criteria to find a trade-off between the time and the value of the objective function.

The last two columns of Table \ref{tbl:Performance} show the time needed for solving \eqref{eq:ReformulatedMaxAll} and \eqref{eq:ReformulatedExp}. The problem \eqref{eq:ReformulatedExp} has more decision variables and requires more time. The timings are insensitive to the stopping criteria of the algorithm but are affected by the stopping criteria of pattern search.  The value of mesh tolerance can be understood as a re-definition of the discretization, as shown in Fig. \ref{fig:MeshDiscrete}. It indicates how close to the current point the pattern search algorithm will look before stopping.  We observe that relaxing the value of mesh tolerance from 0.000001 to 0.01 allows reducing the time to solve \eqref{eq:ReformulatedExp} from around 5 s to 0.1 s and \eqref{eq:ReformulatedMaxAll} from 0.5 s to 0.02.

Conversely, the value of the initial mesh size defines how far from the current point pattern search starts its iterations. We see that increasing the initial mesh size can bring down the number of iterations (here for 10 and the mesh tolerance equal to 0.000001). At the same time, increasing the initial mesh size increases the time to solve \eqref{eq:ReformulatedExp}. This is because pattern search starts looking far from the given initial value. However, as indicated by \cite{Algorithms_Kochenderfer2019}, the expanders are supposed to be at the boundary of the safe set, i.e. close to $y$ where $l_n(y,j)= J_{\min}$ for all $j$. Giving the initial point close to the boundary and using a large initial mesh size requires multiple iterations inside patter search to reduce the mesh size to find a solution. 

As a result, setting the value of the initial mesh size parameter can serve as passing additional information about the function to the pattern search method. The value can also be problem-dependent, as visible from the values obtained for a small initial mesh size and a large one. A large size of the initial mesh may suggest going outside the search area, in particular, if the current point is close to the boundary of $\mathcal{A}$. This phenomenon can be observed when looking at values of $K_{vi}$  obtained for the initial mesh size of 10, and 20. In these cases, the solution for $K_{vi}$ obtained from \eqref{eq:Maxw2} that lies on the boundary. A smaller mesh size, equal to one or five, allowed the solution to lie inside $\mathcal{A}$, resulting in an improved value of the objective.

\subsubsection{Comparison with default SafeOpt}
The performance of the proposed reformulation was then compared to the default version of SafeOpt from Section \ref{sec:DefaultSafeOpt}. We analysed the performance in terms of the time necessary to obtain a solution. The analysis is collected in Table \ref{tbl:ComparisonWithSafeOpt}. The time in the second column is the average time obtained from 10 runs of every algorithm. The stopping criteria used for SafeOpt are the number of iterations and evaluation of all the points in the discretized search space. The parameters of the reformulated algorithm were chosen as $\epsilon_1=\epsilon_2=0.1$, with a mesh tolerance of 0.01 and two initial mesh sizes of 10 (RA 10) and one (RA 1). These values were chosen for the shortest solution time and best value of the objective (bold).

\begin{table}[!tbp]
\centering
\caption{Performance of the default SafeOpt}
\label{tbl:ComparisonWithSafeOpt}
\begin{tabular}{@{}P{0.3cm}cccccP{1.1cm}P{1.1cm}@{}}
\toprule
$N$ & Time [s] & Iter. & $K_p^*$ &$K_v^*$ &$K_{vi}^*$ & \parbox[c]{\hsize}{\centering Obj. Tune} & \parbox[c]{\hsize}{\centering Obj. Sine} \\ \midrule
$10$   &  4 & 3 & 36.7& 11.1& 27.8& -10.2 & -0.14\\
$1000$  & 14  & 4 & 49.7& 29.1& 43.7& -7 & -0.08\\
$2500$  & 56 & 5 &     86.8& 34.3& 41.8     & -4.3 & -0.05\\ \bottomrule
\end{tabular}
\end{table}

The results are collected in Table \ref{tbl:ComparisonWithSafeOpt} and shown in Fig. \ref{fig:SpeedRes1}-\ref{fig:PositionRes1}, corresponding to the trajectory used for tuning, and in Fig. \ref{fig:SpeedResSin1}-\ref{fig:PositionResSin1} for a sinusoidal trajectory, truncated at the position 0.9 cm and zero. In all the cases, the best objective was obtained for the reformulated algorithm RA 1. The default SafeOpt was second best, at the expense of the computational time.

\section{Discussion and conclusions}
\label{sec:Conclusions}
Existing algorithms for safe learning, such as SafeOpt, allow for ensuring safety at the expense of increased computational effort. The current paper proposes a reformulation of the SafeOpt algorithm as a series of optimization problems. By using direct search methods for the optimization problems, we also preserve the derivative-free character of SafeOpt while improving computation time. 

In future work, we plan to analyse the impact of the solver used for optimization problems on the convergence of the algorithm and safety guarantees.

\balance

\bibliography{bibTAC}             % bib file to produce the bibliography

\begin{thebibliography}{12}
\providecommand{\natexlab}[1]{#1}
\providecommand{\url}[1]{\texttt{#1}}
\providecommand{\urlprefix}{URL }
\expandafter\ifx\csname urlstyle\endcsname\relax
  \providecommand{\doi}[1]{doi:\discretionary{}{}{}#1}\else
  \providecommand{\doi}{doi:\discretionary{}{}{}\begingroup
  \urlstyle{rm}\Url}\fi

\bibitem[{{\AA}str{\"o}m and H{\"a}gglund(2006)}]{Advanced_Astroem2006}
{\AA}str{\"o}m, K.J. and H{\"a}gglund, T. (2006).
\newblock \emph{Advanced {PID} Control}.
\newblock ISA-The Instrumentation, Systems, and Automation Society.

\bibitem[{Audet and Hare(2017)}]{Derivative_Audet2017}
Audet, C. and Hare, W. (2017).
\newblock \emph{Derivative-free and blackbox optimization}.
\newblock Springer.

\bibitem[{Berkenkamp et~al.(2021)Berkenkamp, Krause, and
  Schoellig}]{BayesianBerkenkamp2021}
Berkenkamp, F., Krause, A., and Schoellig, A.P. (2021).
\newblock Bayesian optimization with safety constraints: safe and automatic
  parameter tuning in robotics.
\newblock \emph{Machine Learning}, 1--35.

\bibitem[{Berkenkamp et~al.(2016)Berkenkamp, Schoellig, and
  Krause}]{SafeBerkenkamp2016}
Berkenkamp, F., Schoellig, A.P., and Krause, A. (2016).
\newblock Safe controller optimization for quadrotors with {G}aussian
  processes.
\newblock In \emph{2016 {IEEE} International Conference on Robotics and
  Automation ({ICRA})}. {IEEE}.

\bibitem[{Duivenvoorden et~al.(2017)Duivenvoorden, Berkenkamp, Carion, Krause,
  and Schoellig}]{Constrained_Duivenvoorden2017}
Duivenvoorden, R.R.P.R., Berkenkamp, F., Carion, N., Krause, A., and Schoellig,
  A.P. (2017).
\newblock Constrained {B}ayesian optimization with particle swarms for safe
  adaptive controller tuning.
\newblock \emph{{IFAC}-{PapersOnLine}}, 50(1), 11800--11807.

\bibitem[{Fiducioso et~al.(2019)Fiducioso, Curi, Schumacher, Gwerder, and
  Krause}]{Safe_Fiducioso2019}
Fiducioso, M., Curi, S., Schumacher, B., Gwerder, M., and Krause, A. (2019).
\newblock Safe contextual {B}ayesian optimization for sustainable room
  temperature {PID} control tuning.

\bibitem[{Khosravi et~al.(2020)Khosravi, Behrunani, Smith, Rupenyan, and
  Lygeros}]{Cascade_Khosravi2020}
Khosravi, M., Behrunani, V., Smith, R.S., Rupenyan, A., and Lygeros, J. (2020).
\newblock Cascade control: Data-driven tuning approach based on {B}ayesian
  optimization.
\newblock \emph{IFAC-PapersOnLine}, 53(2), 382--387.
\newblock 21th IFAC World Congress.

\bibitem[{Khosravi et~al.(2022)Khosravi, Behrunani, Myszkorowski, Smith,
  Rupenyan, and Lygeros}]{Performance_Khosravi2022}
Khosravi, M., Behrunani, V.N., Myszkorowski, P., Smith, R.S., Rupenyan, A., and
  Lygeros, J. (2022).
\newblock Performance-driven cascade controller tuning with {B}ayesian
  optimization.
\newblock \emph{IEEE Transactions on Industrial Electronics}, 69(1),
  1032--1042.

\bibitem[{Kim et~al.(2021)Kim, Allmendinger, and
  L{\'o}pez-Ib{\'a}{\~{n}}ez}]{Safe_Kim2021}
Kim, Y., Allmendinger, R., and L{\'o}pez-Ib{\'a}{\~{n}}ez, M. (2021).
\newblock Safe learning and optimization techniques: Towards a survey of the
  state of the art.
\newblock In F.~Heintz, M.~Milano, and B.~O'Sullivan (eds.), \emph{Trustworthy
  AI - Integrating Learning, Optimization and Reasoning}, 123--139. Springer
  International Publishing, Cham.

\bibitem[{Kochenderfer and Wheeler(2019)}]{Algorithms_Kochenderfer2019}
Kochenderfer, M.J. and Wheeler, T.A. (2019).
\newblock \emph{Algorithms for optimization}.
\newblock {MIT} Press.

\bibitem[{Lewis and Torczon(2002)}]{Globally_Lewis2002}
Lewis, R.M. and Torczon, V. (2002).
\newblock A globally convergent augmented {L}agrangian pattern search algorithm
  for optimization with general constraints and simple bounds.
\newblock \emph{{SIAM} Journal on Optimization}, 12(4), 1075--1089.

\bibitem[{Sui et~al.(2015)Sui, Gotovos, Burdick, and Krause}]{Safe_Sui2015}
Sui, Y., Gotovos, A., Burdick, J., and Krause, A. (2015).
\newblock Safe exploration for optimization with {G}aussian processes.
\newblock In F.~Bach and D.~Blei (eds.), \emph{Proceedings of the 32nd
  International Conference on Machine Learning}, volume~37 of \emph{Proceedings
  of Machine Learning Research}, 997--1005. PMLR, Lille, France.

\end{thebibliography}
                                                     % with bibtex (preferred)
                                                   
\end{document}